\documentclass[11pt]{amsart}
\usepackage{latexsym, amssymb,lscape,epsfig,setspace}
\theoremstyle{plain}
\newtheorem{theorem}{Theorem}
\newtheorem{proposition}{Proposition}
\newtheorem{conjecture}{Conjecture}

\theoremstyle{definition}

\newtheorem{example}{Example}
\newtheorem{exercise}{Exercise}

\newenvironment{renumerate}%
{%
\begin{enumerate}}%
{\end{enumerate}%
}%

\newenvironment{theo}[1]%
{\begin{theorem}\label{T:#1}}%
{\end{theorem}}

\newenvironment{prop}[1]%
{\begin{proposition}\label{T:#1}}%
{\end{proposition}}

{\begin{conjecture}\label{T:#1}}%
{\end{conjecture}}

\newenvironment{definition}%
{\vskip6pt%
\noindent%
{\bf Definition.}}%
{\vskip6pt}

\newenvironment{remark}%
{\vskip6pt%
\noindent%
{\it Remark.}}%
{\vskip6pt}

{\vskip6pt%
\noindent%
{\it Remarks}. %
\begin{renumerate}}%
{\end{renumerate}\vskip6pt}

{\begin{exercise}\label{T:#1}}%
{\end{exercise}}

\newenvironment{ex}[1]%
{\begin{example}\label{T:#1}}%
{\end{example}}

\newcommand{\RR}{\text{${\mathbb R}$}}
\newcommand{\C}{\text{$\mathbb C$}}
\newcommand{\N}{\text{$\mathbb N$}}
\newcommand{\Z}{\text{$\mathbb Z$}}

\newcommand{\ga}{\text{$\alpha$}}

\newcommand{\gd}{\text{$\delta$}}

\newcommand{\gO}{\text{$\Omega$}}
\newcommand{\go}{\text{$\omega$}}
\newcommand{\gf}{\text{$\varphi$}}

\newcommand{\tensor}{\otimes}
\newcommand{\im}{\mathrm{Im}\,}

\newcommand{\mc}[1]{\text{$\mathcal{#1}$}}

\newcommand{\into}{\rightarrow}
\newcommand{\noqed}{\let\qed\relax}

\begin{document}

\pagestyle{headings}

\title{Formality of $k$-connected spaces in $4k+3$ and $4k+4$ dimensions}
\author{Gil Ramos Cavalcanti}
\address{Mathematical Institute, St. Giles 24-29, Oxford, OX1 3BN, UK}
\email{gilrc@maths.ox.ac.uk}
\thanks{Researcher supported by CAPES (Coordena\c c\~ao de Aperfei\c coamento de Pessoal de N\'ivel Superior, Minist\'erio da Educa\c c\~ao e Cultura), Brazilian Government. Grant 1326/99-6}
%\subjclass{Primary 53D35; Secondary 57R19}
\date{November 14, 2004}
%\keywords{strong Lefschetz property, symplectic blow-up, Massey products}

\begin{abstract}
Using the concept of $s$-formality we are able to extend the bounds of a Theorem of Miller and show that a compact $k$-connected $(4k+3)$- or $(4k+4)$-manifold with $b_{k+1}=1$ is formal. We study $k$ connected $n$-manifolds, $n= 4k+3, 4k+4$, with a hard Lefschetz-like property and prove that in this case if $b_{k+1}=2$, then the manifold is formal, while, in $4k+3$-dimensions, if $b_{k+1}=3$ all Massey products vanish. We finish with examples inspired by symplectic geometry and manifolds with special holonomy.
\end{abstract}

\maketitle
\markboth{\sc G. R. Cavalcanti}{\sc Formality of $k$-connected spaces}

\section*{Introduction}

A theorem of Miller \cite{Mi79} states that any compact orientable $k$-connected manifold of dimension $d \leq (4k+2)$ is formal. In particular, a compact simply-connected $n$-manifold is formal if $n \leq 6$. Recently, Fern\'andez and Mu\~noz gave examples of simply-connected nonformal 7- and 8-manifolds \cite{FM02b} and Dranishnikov and Rudyak gave examples of $k$-connected nonformal manifolds in any dimension greater than $4k+2$ \cite{DrRu03}, therefore proving that Miller's theorem can not be improved without further hypotheses.

Here, we adopt the point of view that for a $k$-connected manifold, the smaller the $k+1^{th}$ Betti number, $b_{k+1}$, is, the simpler the topology. Then we establish the biggest value of $b_{k+1}$ for which one can still assure formality of an $n$-manifold, $n= 4k+3, 4k+4$:

\begin{theo}{miller improved}
A compact orientable $k$-connected manifold of dimension $4k+3$ or $4k+4$ with $b_{k+1} =1$ is formal.
\end{theo}

%The results in this paper arise from my interest on the formality of compact simply-connected 7- and 8-manifolds, which lie just outside the reach of Miller's theorem.

One motivation for the study of such manifolds comes from the existence of special geometric structures in 7- and 8-manifold, e.g., $G_2$, $Spin(7)$ and symplectic structures. A compact irreducible Riemannian manifold with holonomy group $G_2$ or $Spin(7)$ has finite fundamental group and hence its universal cover is compact, simply-connected and has special holonomy. Other examples of manifolds with special holonomy are K\"ahler manifolds which are formal \cite{DGMS75}, hence it is conceivable that there is a connexion between the existence of special holonomy metrics and formality.

Another reason to study formality of 8-manifolds comes from symplectic geometry, as this is the only dimension where the question of existence of compact 1-connected nonformal examples is still open.

One property shared by irreducible $G_2$- and $Spin(7)$-manifolds is that they have a hard Lefschetz-like property. If $M^n$ is one such manifold, there is a closed $(n-4)$-form \gf\ for which:
\begin{equation}\label{E:motivation}
[\gf] \cup :H^2(M,\RR) \stackrel{\cong}{\into} H^{n-2}(M,\RR),
\end{equation}
is an isomorphism \cite{Jo00}. Of course, this last property is also shared by 2-Lefschetz symplectic 8-manifolds, where, by definition, \eqref{E:motivation} holds with $\gf = \go^2$.

We are interested in how \eqref{E:motivation} can be used to improve on Theorem \ref{T:miller improved}. So we prove that if $M^n$, $n =4k+3,4k+4$, is $k$-connected and there is $\gf \in H^{n-2k-2}(M)$ for which
\begin{equation}\label{E:gf and h2}
\gf \cup :H^{k+1}(M,\RR) \stackrel{\cong}{\into} H^{n-k-1}(M,\RR)
\end{equation}
is an isomorphism and $b_{k+1}(M)=2$, then, $M$ is formal. If $n= 4k+3$ and $b_{k+1}(M)=3$, then all Massey products on $M$ vanish uniformly. We construct examples showing that our bounds are optimal.

This paper is organized as follows. In the first section we cover background material: introduce minimal models, formality, $s$-formality and Massey products. The proof that every compact 1-connected manifold has a minimal model is given to fix notation. In Section \ref{miller}, we prove Theorem \ref{T:miller improved} by explicitly constructing a minimal model. In Section \ref{lefschetz}, we prove that the existence of the isomorphism \eqref{E:gf and h2} has a `formalizing tendency' in a $k$-connected orientable $4k+3$- or $4k+4$-manifold, as it implies formality for $b_2 =2$ and, in the $4k+3$-dimensional case, vanishing of Massey products for $b_2 =3$. In the last section we study examples from symplectic and Riemannian geometry where our results can be applied.

I would like to thank Nigel Hitchin for his guidance and Marisa Fern\'andez and Yuli Rudyak for their suggestions and encouragement.

\section{Minimal models and $s$-formality}\label{s-formality}

In this section we introduce the basic concepts we are interested in: differential graded algebras, minimal models, formality and $s$-formality. All the results are standard (see \cite{GM81,FHT01}). Some proofs are given in order to fix notation.

\begin{definition}
A {\it differential graded algebra}, or {\it DGA} for short, is an $\N$ graded vector space $\mc{A}^{\bullet}$ over a field $k$, endowed with a product and a differential $d$ satisfying:
\begin{enumerate}
\item The product maps $\mc{A}^{i}\times \mc{A}^{j}$ to $\mc{A}^{i+j}$ and is graded commutative:
$$a \cdot b = (-1)^{ij} b\cdot a;$$
\item The differential has degree 1: $d:\mc{A}^{i} \into \mc{A}^{i+1}$;
\item $d^2 =0$;
\item The differential is a derivation: for $a \in \mc{A}^{i}$ and $b \in \mc{A}^{j}$
$$d(a\cdot b) = da \cdot b + (-1)^i a \cdot db.$$
\end{enumerate}
\end{definition}

We will deal only with DGAs over \RR. A nontrivial example of DGA is the complex of differential forms on a manifold equipped with the exterior derivative. Another example is the cohomology of a manifold, with $d=0$ as differential. The cohomology of a DGA is defined in the standard way.

%Sometimes, given a DGA, $(\mc{A},d)$, one can construct another differential  graded algebra that captures all the information about the differential and which is minimal in the following sense.
%If the fundamental group of this manifold is nilpotent, then this algebra contains all the information about the rational homotopy type of the manifold: $\pi_i(M)\tensor \Q$, see \cite{Su78}.

\begin{definition}
A DGA $(\mc{M},d)$ is {\it minimal} if it is free as a DGA (i.e. polynomial in even degree and skew symmetric in odd degree) and has generators $e_1, e_2 \dots, e_n, \dots$ such that:
\begin{enumerate}
\item \label{formalitypage}The degree of the generators form a weakly increasing sequence of positive numbers;
\item There are finitely many generators in each degree;
\item The differential satisfies $de_i \in \wedge\{e_1, \dots, e_{i-1}\}$.
\end{enumerate}
A {\it minimal model} for a differential graded algebra $(\mc{A},d)$ is a minimal DGA, $\mc{M}$, together with a map of differential graded algebras, $\rho: \mc{M} \into \mc{A}$, inducing isomorphism in cohomology. If $(\mc{A},d)$ is the complex of differential forms on a manifold $M$ we also refer to the minimal model of $\mc{A}$ as the minimal model of $M$.
\end{definition}

Every simply-connected compact manifold has a minimal model which is unique up to isomorphism. As the cohomology algebra of a manifold with real coefficients is also a DGA we can also construct its minimal model. This minimal model is not necessarily the same as the minimal model for the manifold and in general it carries less information.

\begin{definition}
A manifold is {\it formal} if the minimal models for $M$ and for its cohomology algebra are isomorphic, or equivalently, if there is a map of differential graded algebras $\psi:\mc{M} \into H^{\bullet}(M,\RR)$ inducing isomorphism in cohomology, where $\mc{M}$ is the minimal model of $M$.
\end{definition}

The proof that every 1-connected manifold has a minimal model is done in a constructive way and the fundamental tool is a Hirsch extension. If $\mc{A}$ is a differential algebra,
$V$ a finite dimensional vector space and $d:V \into
\ker(d)\cap (\mc{A}^{k+1})$ a linear map, we can form the {\it Hirsch extension}
$\mc{A}\tensor_d (\wedge V)_k$ as follows: the elements of
$\mc{A}\tensor_d (\wedge V)_k$ are elements of the tensor product of
$\mc{A}$ with the free algebra generated by $V$ in degree $k$, the
differential is defined by linearity plus the product rule
$d(a\tensor v) = da \tensor v + (-1)^{\deg(a)} a \wedge dv \tensor
1$. We denote 
$1\tensor v$ by $v$ and $a \tensor 1$ by $a$ and then $a \tensor v$
becomes simply $a \wedge v$.

In doing such an extension one is, amongst other things,

\begin{renumerate}
\item killing cohomology in degree $k+1$ since the elements in the image
of $d$  that were closed in $\mc{A}$, but not exact, become exact in $\mc{A}\tensor_d (\wedge V)_k$;

\item creating cohomology in degree $k$, namely, the new classes are
given by the kernel of $d:V \into \mc{A}^{k+1}/\im(d:\mc{A}^{k} \into \mc{A}^{k+1})$.
\end{renumerate}

Therefore, we can always split $V = B^k \oplus Z^k$, where the map $d:B^k \into \mc{A}^{k+1}/\im(d:\mc{A}^{k} \into \mc{A}^{k+1})$ is an injection and $Z_k = \ker(d:V \into \mc{A}^{k+1})/\im(d:\mc{A}^{k} \into \mc{A}^{k+1})$, although this splitting is not canonical.

\begin{theo}{existence of minimal models}{\em (Sullivan \cite{Su78})}
Every compact simply-connected manifold has a minimal model.
\end{theo}
\begin{proof}
The proof is accomplished by constructing the minimal model, which is done by induction. The starting point is a free DGA, $\mc{M}_2$, generated by $H^2(M)$ in degree 2 and with vanishing differential. We define $\rho_2: \mc{M}^{\bullet}_2 \into \gO^{\bullet}(M)$ by choosing arbitrary representatives for the cohomology classes in $H^2(M)$ and extend it to higher symmetric powers of $H^2(M)$ so that it is a map of algebras.

Now assume we have a free DGA, $\mc{M}_k$, generated by elements of degree at most $k$  with a map of DGAs $\rho_k :\mc{M}^{\bullet}_k \into \gO^{\bullet}(M)$ such that 
$\rho_k^*:H^i(\mc{M}_k) \into H^i(M)$ is an isomorphism for $i \leq k$ and an injection for $i=k+1$. Let 
$$V^{k+1} = \ker(\rho_k^*: H^{k+2}(\mc{M}_k) \into H^{k+2}(M)) \oplus H^{k+1}(M)/\rho_k^*( H^{k+1}(\mc{M}_k)).$$
Choosing linearly a representative $v$ for each cohomology class
$$[v] \in \ker(\rho_k^*: H^{k+2}(\mc{M}_k) \into H^{k+2}(M)),$$
define $\mc{M}_{k+1}$ to be the Hirsch extension of $\mc{M}_{k}$ by $V^{k+1}$, with $d$ given by:
$$\begin{cases}
d[v] = v & \mbox{ if } [v] \in \ker(\rho_k^*: H^{k+2}(\mc{M}_k) \into H^{k+2}(M),\\
dv = 0 & \mbox{ if } v \in H^{k+1}(M)/\rho_k^*( H^{k+1}(\mc{M}_k)).
\end{cases}$$
Finally, we let $\rho_{k+1}|_{\mc{M}_k} =\rho_k$ and define it in $V^{k+1}$ in the following way. For $[v] \in \ker(\rho_k^*: H^{k+2}(\mc{M}_k) \into H^{k+2}(M))$, by definition, $\rho_{k}(d[v]) = da$, for some $k+1$-form $a$. We let $\rho_{k+1}([v])= a$. For $v \in  H^{k+1}(M)/\rho_k^*(H^{k+1}(\mc{M}_k))$, we choose a closed form $b$ whose cohomology class represents the same element and let $\rho_{k+1}(v) =b$. These choices can be made so the $\rho_{k+1}:V^{k+1} \into \gO^{k+1}(M)$ is linear. The map $\rho_{k+1}$ is extended to the rest of $\mc{M}_{k+1}$ by requiring it to be a map of DGAs.

With these choices, $\mc{M}_{k+1}$ is a free DGA generated by elements of degree at most $k+1$ with a map of DGAs $\rho_{k+1} :\mc{M}_k \into \gO^{\bullet}(M)$ inducing an isomorphism $\rho_{k+1}^*:H^i(\mc{M}_{k+1}) \into H^i(M)$ for $i \leq k+1$ and an injection for $i=k+2$ and the inductive step is complete.

This procedure gives us a family of free DGAs, each obtained by a Hirsch extension from the previous one:
$$ \mc{M}_2 < \mc{M}_3 < \cdots \mc{M}_k < \cdots,$$
and maps $\rho_k:\mc{M}_k \into \gO^{\bullet}(M)$ such that:
\begin{renumerate}
\item $\rho_k|{\mc{M}_i} = \rho_i$ for $i\leq k$, hence induce a map in the limit DGA $\rho: \mc{M} \into \gO^{\bullet}(M)$;
\item The induced maps in cohomology are isomophisms $\rho^*: H^i(\mc{M}_k) \cong H^i(M)$, for $i \leq k$ and therefore, in the limit, $H^i(\mc{M}) \cong H^i(M)$ for all $i$.
\end{renumerate}
Thus the limit $\mc{M}$ of the $\mc{M}_k$ is the minimal model for $M$.
\end{proof}

The main reason to present the proof of this theorem is that it is a constructive one and sometimes properties about the minimal model are more easily obtained using this point of view. For example we have an alternative description of formality (see \cite{FHT01,GM81}, for a proof).

\begin{prop}{alternativeformal}
Let $V^k$ be the spaces introduced in degree $k$ when constructing the minimal model of a manifold $M$. $M$ is formal if and only if there is a splitting $V^k = B^k \oplus Z^k$, such that
\begin{enumerate}
\item $d:B^k \into \mc{M}^{k+1}$ is an injection;
\item $d(Z^k) = 0$;
\item If an element of the ideal $\mc{I}(\oplus_{k \in \N} B^k) < \mc{M}$ is closed, then it is exact.
\end{enumerate}
\end{prop}

This characterization of formality allows one to consider weaker versions. Notably, Fern\'andez and Mu\~noz introduced in \cite{FM02} the useful concept of $s$-formality.

\begin{definition}
A manifold is $s$-formal if there is a choice of splitting $V^k = Z^k \oplus B^k$ satisfying (1) and (2) above and such that every closed element in the ideal $\mc{I}(\oplus_{k \leq s}B^k) < \mc{M}_s$ is exact in $\mc{M}$.
\end{definition}

Clearly this is a weaker concept, in general, but it is also obvious that if an $n$-dimensional manifold is $n$-formal, it is formal. The following result of Fern\'andez and Mu\~noz shows that sometimes the weaker condition of $s$-formality implies formality.

\begin{theo}{fernandezmunoz}{\em (Fern\'andez and Mu\~noz \cite{FM02})}: If a compact orientable manifold $M^n$ is $s$-formal for $s \geq n/2 -1$, then $M$ is formal.
\end{theo}

The point of this theorem is that in order to prove formality by constructing the minimal model and finding the splitting as in Proposition \ref{T:alternativeformal}, it is only necessary to determine the beginning of the minimal model.

So far we have two ways to tell whether a manifold is formal: from the definition of formality and from a (partial) construction of the minimal model and the theorem of Fern\'andez and Mu\~noz. Next we introduce Massey products, which are an obstruction to formality. The ingredients are  $a_{12},a_{23},a_{34} \in
\mc{A}$ three closed elements such that $a_{12}a_{23}$ and $a_{23}a_{34}$ are exact. Then, denoting $\bar{a}= (-1)^{|a|}a$, we define
\begin{equation}\tag{$*$}
\begin{cases}
\overline{a_{12}}a_{23} &= da_{13}\\
\overline{a_{23}}a_{34} &= da_{24}.
\end{cases}
\end{equation}
In this case, one can consider the element
$\overline{a_{13}} a_{34} + \overline{a_{12}} a_{24}$. By the choice of $a_{13}$ and $a_{24}$ this element is closed, hence it represents a cohomology class. Observe,
however, that $a_{13}$ and $a_{24}$ are not well defined and we can change
them by any closed element, hence the expression above does not define a
unique cohomology class but instead an element in the
quotient $H^{\bullet}(\mc{A})/\mathcal{I}([a_1],[a_3])$.

\begin{definition}
The {\it Massey product} or {\it triple Massey product}
$\langle [a_{12}], [a_{23}], [a_{34}]\rangle$, 
of the cohomology classes $[a_{12}]$, $[a_{23}]$ and $[a_{34}]$
with $[a_{12}] [a_{23}]=[a_{23}] [a_{34}]=0$
 is the coset
$$\langle [a_{12}],[a_{23}],[a_{34}]\rangle = [\overline{a_{12}}
  a_{24} + \overline{a_{13}} a_{34}] + ([a_{12}], [a_{34}])
  \in H^{\bullet}(\mc{A})/\mathcal{I}([a_{12}],[a_{34}]),$$
where $a_{13}$ and $a_{24}$ are defined by $(*)$. 
We say that all the Massey products {\it vanish uniformily} if there is a fixed set of choices for which all possible Massey products are represented by exact forms.
\end{definition}

The importance of Massey products for this work comes from the following result: {\it 
If $M$ is formal, all the Massey products vanish uniformily.}

\iffalse
\begin{renumerate}
\item If $\rho:\mc{M} \into \mc{A}$ is the minimal model for $\mc{A}$ and we have a Massey product $\langle v_{12}, \cdots, v_{n-1n} \rangle$ in $\mc{A}$, we can use the same cohomology classes to get a Massey product in $\mc{M}$. As $\rho$ induces an isomorphism in cohomology, the Massey product vanishes in $\mc{A}$ if and only if it vanishes in $\mc{M}$. In some sense, the minimal model is the natural place in which to define Massey products;
\item Using the notation of Therorem \ref{T:alternativeformal}, there is a decomposition of $\mc{M}$, as a vector space, into
$$ \mc{M} = \mbox{span}(\oplus Z_i) \oplus \mc{I}(\oplus B_i);$$
where every element in the first summand is closed. If an element $a$ is exact, say $a= db$, we can split $b$ according to the decomposition above into $b_1 + b_2$, but the element $b_1 \in  \mbox{span}(\oplus Z_i)$ is closed, so in fact we have that every exact element satisfies:
$$a = d b_2 \in d(\mc{I}(\oplus B_i)).$$
\end{renumerate}

Now assume that $M$ is formal. Then any Massey product gives rise to a product in the minimal model from ({\it i\,}) and we can always make choices so that it lies in $\mc{I}(\oplus B_i)$, from ({\it ii\,}), but, by formality, this implies that the Massey product vanishes.

\begin{prop}{vanishing massey products and formality}
If a compact manifold has nonvanishing Massey products, it is not formal.
\end{prop}
\fi

\section{Extending Miller's bounds}\label{miller}

In this section we prove Theorem \ref{T:miller improved}. 
% and give examples of nonformal simply-connected 7- and 8-manifolds with $b_2=2$, showing that the bounds given by that theorem are sharp in these dimensions.
The way to prove formality in this case is by constructing the beginning of the minimal model for the manifold and then using the theorem of Fern\'andez and Mu\~noz (Theorem \ref{T:fernandezmunoz} above).

\vskip6pt
\noindent
{\sc Proof of Theorem \ref{T:miller improved}}. We will only prove the theorem for $4k+4$-manifolds, as the other case is analogous. The proof is accomplished by constructing the minimal model.

As $M$ is $k$-connected, $\mc{M}_{2k}$ is a free algebra generated by $H^{i}(M)$ in degree $i \leq 2k$ with vanishing differential and $\rho$ maps linearly each cohomology class to a representative form. The first time we may have to use a nonzero differential (and hence introduce one of the $B^j$ spaces) is in degree $2k+1$. This will be the case only if the generator $a \in \mc{M}_{k+1}$ satisfies $a^2 \neq 0$ and $\rho(a^2)$ is an exact form. Hence if either $\rho(a^2)$ is not exact or $a^2=0$,  all the spaces $B^j$ are trivial for $j \leq 2k+1$, showing that $M$ is $2k+1$-formal and hence, by Theorem \ref{T:fernandezmunoz}, formal.

So we only have to consider the case where the cohomology class $\ga \in H^{k+1}(M)$ satisfies $\ga^2 =0$ and $k+1$ is even. In this case, the Hirsch extension in degree $2k+1$ is given by
$$V^{2k+1} = H^{2k+1}(M)\oplus \mbox{span}\{b\},$$
where $d$ vanishes in $H^{2k+1}(M)$, $db = a^2$ and $\rho(b)$ is a form such that $d\rho(b) = \rho(a^2)$. With this splitting, $\mc{I}(\oplus_{j \leq 2k+1}(B^j)) < \mc{M}_{2k+1}$ is just the ideal generated by $b$ and to prove formality using Theorem \ref{T:fernandezmunoz} we have to show that any closed form in this ideal is exact in $\mc{M}$.

A closed form in this ideal, being the product of $b$ and an element of degree at least $k+1$, will have degree at least $3k+2$. Since $M$ is $k$-connected, Poincar\'e duality gives $H^{j}(M)=\{0\}$ for $3k+3 < j < 4k+4$. If an element in $\mc{I}(b)$ of degree $3k+2$ or $3k+3$ is closed and nonexact in $\mc{M}$, Poincar\'e duality implies that its dual is in either $\mc{M}_{k+2}^{k+2} \cong H^{k+2}(M)$ (in the former case) or in $\mc{M}_{k+1}^{k+1} \cong H^{k+1}(M)$ (in the latter). Either way, from there we can produce a degree $4k+4$ closed nonexact element in $\mc{I}(b)$. So we only have to check that any $4k+4$ closed element in $\mc{I}(b)$ is exact.

The only elements in $\mc{I}(b)$ in degree $4k+4$ are of the form $a b v$, with $v \in H^{k+2}(M) \cong V^{k+2}$, which have derivative $d(a b v) = a^3v$ and hence are not closed if $v \neq 0$. This shows that $M$ is $(2k+1)$-formal and therefore formal.
\qed
\vskip6pt

\section{Formality of hard Lefschetz manifolds}\label{lefschetz}

In this section we study compact orientable $k$-connected $n$-manifolds, $n=4k+3, 4k+4$, for which there is a cohomology class $\gf \in H^{n-2k-2}(M)$ inducing an isomorphism
\begin{equation}\tag{\ref{E:gf and h2}}
\gf \cup :H^{k+1}(M) \stackrel{\cong}{\into} H^{n-k-1}(M).
\end{equation}
We prove that this property has a `formalizing tendency' in the following sense. 

\begin{theo}{g2 with b2 leq 3} 
A compact orientable $k$-connected manifold $M^n$, $n=4k+3, 4k+4$, satisfying \eqref{E:gf and h2} with $b_{k+1}=2$ is formal. If $n=4k+3$ and $b_{k+1}=3$ all the Massey products vanish uniformly.
\end{theo}
\begin{proof}
We start considering the case $b_{k+1}=2$. The cases $n=4k+3$ and $n=4k+4$ are similar, so we only deal with the latter. Due to \eqref{E:gf and h2}, the class $\gf$ induces a nondegenerate bilinear form on $H^{k+1}(M)$. If $k$ is even, this bilinear form is skew, and if $k$ is odd, it is symmetric.

For $k$ odd, the signature of the bilinear form, i.e., the difference between the number of positive and negative eigenvalues, is either $2$, $0$ or $-2$. By changing $\gf$ to $-\gf$, the case of signature $-2$ can be transformed into signature 2, hence there are two possibilities to consider. As they are similar, we will only treat the signature 2 case.

Let $a_1,a_2$ ($da_1=da_2=0$) be generators of $\mathcal{M}_{k+1} = Sym^{\bullet}H^{k+1}(M)$, the
 first nontrivial stage of the minimal model for $M$. We may further assume that the bilinear form induced by \gf\ is diagonal in the basis $\{a_1,a_2\}$ and
$$ \int \gf a_i^2 =1.$$
As with Theorem \ref{T:miller improved}, $\mc{M}^{2k}$ is just the free algebra generated by $H^j(M)$ in degree $j \leq k$ with vanishing differential and the first time we may have to introduce one of the $B^j$ spaces is in degree $2k+1$, where
$$B^{2k+1} \cong \ker Sym^2H^{k+1}(M) \stackrel{\cup}{\into} H^{2k+2}(M).$$

We know that both
 $a_1^2$ and $a_2^2$ are nonzero in $H^{2k+2}(M)$, since by \eqref{E:gf and h2} they pair nontrivially with \gf. So
$$\dim\,B^{2k+1} = \dim(\ker Sym^2H^{k+1}(M) \stackrel{\cup}{\into} H^{2k+2}(M)) \leq \dim(Sym^2 H^{k+1}(M)) - 1 = 2,$$
and hence there may be at most
 two generators in degree $2k+1$ in $\mathcal{M}_{2k+1}$ to kill cohomology
 classes in $\mathcal{M}_{2k}$ that are not present in $H^{2k+2}(M)$. The case
 when only one generator is added has a proof very similar to the one
 of Theorem \ref{T:miller improved}, so we move on to the case when there are two generators $b_1$ and $b_2$ added to kill cohomology in degree $2k+2$.

Then
$$db_i = k^i_{11}a_1^2 + k^i_{12}a_1 a_2 + k^i_{22}a_2^2,\qquad i=1,2.$$
Multiplying by \gf\ and integrating we get $k^i_{11} =
-k^i_{22}$. Hence, by re-scaling and taking linear combinations we may assume that $db_1= a_1^2 - a_2^2$ and $db_2 = a_1 a_2$.

Now, if $c \in \mc{I}(b_1, b_2)_{\leq 2k+1}$ is a closed element (again, by Poincar\'e duality we may assume it has
degree $n$) we can write it as
$$ c = (a_1 c_{11} + a_2c_{12}) b_1 + (a_1 c_{21} + a_2c_{22}) b_2 + (k_1a_1 + k_2 a_2)b_1b_2,$$
where $c_{ij}$ are closed elements of degree $k+2$ and $k_i \in \RR$. The condition $dc=0$ implies that $c=0$, thus every closed form in $(b_1,
b_2)_{\leq 2k+1}$ is exact and $M$ is $2k+1$-formal and therefore formal.

The case when $k$ is even is easier. Indeed, the argument above can be used again, but with $Sym^2 H^{k+1}(M)$ replaced by $\wedge^2 H^{k+1}(M)$ and hence the nondegeneracy of the pairing implies
$$\dim\,B^{2k+1} = \dim(\ker \wedge^2 H^{k+1}(M) \stackrel{\cup}{\into} H^{2k+2}(M)) \leq \dim(\wedge^2 H^{k+1}(M)) - 1 = 0,$$
hence $M$ is trivially $2k+1$-formal and therefore formal.
\vskip6pt

To finish the proof, we consider a $4k+3$-manifold with $b_{k+1}=3$ and prove that all triple Massey products vanish if \eqref{E:gf and h2} holds. Initially we remark that $b_{k+1}=3$ and \eqref{E:gf and h2} can not happen if $k$ is even, as there is no nondegenerate skew bilinear form in an odd dimensional vector space.

We also observe that the Massey product $\langle a,b,c \rangle$ has
degree at least $3k+2$ and since $H^j(M)=\{0\}$, for $3k+2 < j < 4k+3$, this
product will vanish, whenever defined, if its degree is neither $3k+2$ nor $4k+3$. If $\langle a,b,c \rangle \in H^{4k+3}(M)$, it lies in the ideal
generated by $([a],[c])$, so the product also vanishes, therefore the only case left is when $a$, $b$ and $c$ have degree $k+1$ and the product lies in $H^{3k+2}(M)$.

This product vanishes trivially if $c= \lambda a$, hence we can assume $a$ and $c$ linearly independent. Since $\gf a b$ and $\gf b c$ vanish, but \gf\ induces a nondegenerate bilinear form, there is $\beta \in H^{k+1}(M)$ such that $\gf b \beta \neq 0$ and we can express the Massey product in the basis $\{\gf [a],\gf [\beta],\gf [c]\}$:
$$\langle a,b,c \rangle = -v_1 c + a v_2= k_1 \gf[a] + k_2\gf[\beta]+k_3\gf[c].$$
where $ab = dv_1$ and $bc= dv_2$. Multiplying the equation above by $-[b]$ and integrating over $M$ we get
$$ - k_2 \int \gf b \beta  = \int  - v_1 c b + a v_2 b = \int - v_1dv_2 + v_2dv_1 = \int d(v_1v_2) =0.$$
Thus $k_2=0$ and $\langle a,b,c \rangle = k_1 \gf[a] + k_3\gf[c] \in ([a],[c])$. So the Massey product vanishes.

Observe that if $ac = dv_3$, it is still possible to choose $v_i$, $i=1,2,3$ in such a way that the Massey products $\langle a,b,c \rangle$, $\langle b,c,a \rangle$ and $\langle c,a,b  \rangle$ vanish simultaneously. Indeed, let us assume $a$, $b$ and $c$ are linearly independent (the linearly dependent case is similar) and a set of choices of $v_i$ was made:
$$
\begin{cases}
\langle a,b,c \rangle = [v_1 c - a v_2] = k_1 \gf[a] + k_2\gf[c],&\\
\langle b,c,a \rangle = [v_2 a - b v_3] = l_1 \gf[b] + l_2\gf[a],&\\
\langle c,a,b \rangle = [v_3 b - c v_1] = m_1 \gf[c] + m_2\gf[b].&
\end{cases}
$$
Then we can set $v_1 = \tilde{v}_1 + k_2 \gf$, $v_2 = \tilde{v}_2 + l_2 \gf$ and $v_3 = \tilde{v}_3+ m_2 \gf$. With these choices we get
$$
\begin{cases}
\langle a,b,c \rangle = [\tilde{v}_1 c - a \tilde{v}_2] = \tilde{k} \gf[a],&\\
\langle b,c,a \rangle = [\tilde{v}_2 a - b \tilde{v}_3] = \tilde{l} \gf[b],&\\
\langle c,a,b \rangle = [\tilde{v}_3 b - c \tilde{v}_1] = \tilde{m} \gf[c].&
\end{cases}
$$
Adding them up, the left hand side vanishes, giving
$$0 = \tilde{k}\gf [a] + \tilde{l} \gf [b] + \tilde{m} \gf [c].$$
Since $\gf[a]$, $\gf[b]$ and $\gf[c]$ are linearly independent,
$\tilde{k}, \tilde{l}$ and $\tilde{m}$ vanish and, with these choices,
all the Massey products vanish {\it simultaneously}, hinting
at formality.
\end{proof}

\section{Examples}

In this section give simple examples showing that the bounds established in Theorem \ref{T:miller improved} are sharp. We also apply our results to the blow-up of $\C P^8$ along a symplectic submanifold and to Kovalev's examples of $G_2$-manifolds. We finish with an example of a compact 1-connected 7-manifold which satisfies all known topological restrictions imposed by a $G_2$ structure, has $b_2 =4$ and nonvanishing Massey products. This last example shows that the results of Theorem \ref{T:g2 with b2 leq 3} are sharp in 7 dimensions and that one can not answer the question of formality of $G_2$-manifolds using only the currently known topological properties of those.

\begin{ex}{sharp}
In this example we show that the bounds obtained in Theorem \ref{T:miller improved} are sharp in two ways: we construct nonformal compact $k$-connected manifolds $M^n$ with $b_{k+1}=2$, for $n = 4k+3, 4k+4$ and  with $b_{k+1}=1$ for $n > 4k+4$. 

We shall start dealing with odd dimensional examples. Let $k_1 \leq  k_2 \in \N$ be two natural numbers such that $k_1+k_2$ is even and let $V$ be a $(k_1+k_2+2)$-dimensional vector bundle over $S^{k_1+1} \times S^{k_2+1}$ whose Euler class $\chi$ is a nonvanishing top degree cohomology class. Let $X^{2k_1+2k_2+3}$ be the total space of the sphere bundle of $V$. We claim that $X$ has a nontrivial Massey product. Indeed, using the Gysin sequence for this sphere bundle
$$ \cdots \stackrel{\chi\cup}{\into} H^j(S^{k_1+1}\times S^{k_2+1}) \stackrel{\pi^*}{\into} H^j(X) \stackrel{\pi_*}{\into}  H^{j-k_1-k_2-1}(S^{k_1+1}\times S^{k_2+1}) \stackrel{\chi\cup}{\into} H^{j+1}(S^{k_1+1}\times S^{k_2+1}) \into \cdots$$
we see that $H^{k_1+1}(X) = \mbox{span}\{v_1\}$, $H^{k_2+1}(X) =\mbox{span}\{v_2\}$ and $H^{k_1+k_2+2}(X) = \{0\}$, where $v_i$ are generators for the top degree cohomology of each sphere. Therefore, $v_1 \cup v_2 =0$. If $\go_i$ are volume forms pulled back from each $S^{k_i+1}$, then $\go_1 \wedge \go_2 = d\xi$, where
$$\int_{S^{k_1+k_2+1}}\xi = \int_{S^{k_1+1} \times S^{k_2+1}}\chi.$$
Therefore we can compute the Massey product
$$ \langle v_1, v_2, v_2 \rangle = - [\xi \wedge \go_2].$$
This is not exact, as it pairs nontrivially with $v_1$,
$$\int_X \xi \wedge \go_1 \wedge \go_2 =  \int_{S^{k_1+1} \times S^{k_2+1}} \chi \cdot \int_{S^{k_1+1}} \go_1 \int_{S^{k_2+1}}\go_2  \neq 0.$$
And this Massey product has no indeterminacy, hence $X$ is not formal.

It is clear that $X$ is $k_1$-connected and has dimension $2k_1 + 2k_2 +3 \geq 4k_1 +3$. If $k_1 < k_2$, then $b_{k_1+1}(X)=1$, and if $k_1= k_2=k$, then $b_{k+1}(X)=2$.

To obtain even dimensional examples we consider $k_1$ and $k_2$ as before, but now we take a vector bundle $V$ of rank $k_1+k_2+2$ over the manifold $(S^{k_1+1}\times S^{k_2+1}) \# \C P^{k+1}$, where $\#$ indicates connected sum and $k = \frac{1}{2}(k_1+k_2)$. We can further assume that $V$ is trivial outside a small disc. As before, the total space, $X$, of the sphere bundle of $V$ yields a nonformal $(2k_1 + 2k_2+3)$-manifold. The total space of a principal circle bundle over $X$ with Chern class equal to the generator of the second cohomology  $H^{2}(\C P^{k+1},\Z)$ will be a nonformal $k_1$-connected $(2k_1 + 2k_2+4)$-manifold.
\end{ex}

\begin{ex}{symplectic}
Blowing up $\C P^n$ along a suitable submanifold, Babenko and Taimanov proved  in \cite{BT00b} that there are compact 1-connected nonformal symplectic manifolds in any dimension $2k \geq 10$. Due to Miller's theorem, such examples do not exist in dimensions 6 or less and the question is still open for 8-manifolds.

Using the techniques of \cite{Ca04}, one can show that the blow-up, $X$, of $\C P^n$ along any submanifold is always 2-Lefschetz and, if the submanifold is connected, $b_2(X) = 2$. Therefore, our results imply  the blow-up of $\C P^8$ along a connected symplectic submanifold is always formal. One can also try to blow up $\C P^n$ and then take a sequence of Donaldson submanifolds until the result is an 8-manifold, but the manifold obtained this way is also 2-Lefschetz \cite{FM02} and has $b_2 =2$ \cite{Do96}.
\end{ex}

\begin{ex}{kovalev}
In \cite{Kov00}, Kovalev produces a series of examples of 1-connected compact $G_2$-manifolds from pairs of Fano 3-folds via twisted connected sums. If $M^7$ is obtained from the Fanos $F_1$ and $F_2$, he proves that $b_2(M) \leq \min\{b_2(F_1),b_2(F_2)\} -1$. Our results imply that if either of the Fanos involved have $b_2(F_i) \leq 3$, $M^7$ is formal, while if $b_2(F_1)=4$ and $b_2(F_2) \geq 4$, the Massey products vanish. Hence the only possibility for one of his examples to have nontrivial Massey products is if it is constructed from 2 Fanos with $b_2(F_i) \geq 5$.

According to the classification of Fano 3-folds \cite{MoMu86}, these have $b_2 \leq 10$ and if $b_2 \geq 6$, the Fano is just the blow-up of $\C P^3$ in an appropriate number of points. It is easy to follow Kovalev's construction to prove that if one of the summands is $\C P^3$ with some points blown-up, $M^7$ is formal. The case where both the Fanos have $b_2=5$ is more difficult, but one still has formality.
\end{ex}

In \cite{Jo00}, Joyce shows that a compact Riemannian manifold with holonomy $G_2$ has finite $\pi_1$, nonvanishing first Pontryagin class $p_1$ and, if $\gf$ is the closed 3-form determining the structure,
\begin{equation}\label{E:g2 conditions}
%\begin{aligned}\label{E:g2 conditions}
\int a^2 \gf < 0 \mbox{ for } a \in H^{2}(M)\backslash \{0\} \qquad \mbox{ and } \qquad \int p_1 \wedge [\gf ] < 0.
%\end{aligned}
\end{equation}

Using circle bundles we can construct a nonformal manifold satisfying all the topological properties above. The key is Wall's classification of 1-connected spin 6-manifolds.

\begin{theo}{Wall}
{\em (Wall \cite{Wa66})} Diffeomorphism classes of oriented 6-manifolds with
torsion-free homology and vanishing second Stiefel-Witney class
correspond bijectively to isomorphism classes of systems of
invariants:
\begin{itemize}
\item Two finitely generated free abelian groups $H^2$, $H^3$, the latter of even rank;
\item A symmetric trilinear map $\mu:H^2\times H^2\times H^2 \into Z$;
\item A homomorphism $p_1:H^2 \into Z$;
\item Subject to: for $x, y \in H^2$,
$$\mu(x,x,y) = \mu(x,y,y) ~~~\mod(2),$$
\noindent
\phantom{Subject to:} for $x \in H^2$,
$$p_1(x) = 4 \mu(x,x,x) ~~~ \mod(24). $$ 
\end{itemize}
\end{theo}

With appropriate choices for the pairing $\mu$ and for the Chern class of the principal circle bundle, we can obtain nonformal 7-manifolds satisfying \eqref{E:g2 conditions}. As the base manifold is spin, so will be the total space of the circle bundle. We finish this paper with one example constructed this way.

\begin{ex}{second 7manifold}
We let $H^2 = \langle\go, \ga_1,\ga_2,\ga_3,\gamma\rangle$ and
define the cup product on $H^2$ so as to have the following relations
\begin{center}
\vskip6pt
\noindent
\begin{tabular}{c c c c c}
$\go \gamma =0,$& $\go \ga_i \ga_j
= 2\gd_{ij},$& $\go^2 \ga_i=0,$& $\go^3=2,$ & $\go \ga_1 = \gamma \ga_3,$\\
&$\ga_1\ga_2=0,$& $\ga_2 \ga_3 =\gamma\ga_1$& and& $\ga_3 \go = \gamma \ga_1$.
\end{tabular}
\vskip6pt
\end{center}
One set of choices that gives the desired result is the following
\vskip6pt
\noindent
\begin{alignat*}{6}
\go \ga_i \ga_j &= 2\delta_{ij}&\qquad \ga_1\ga_2^2 &=0 &\qquad \ga_1^2 \ga_2 &= 0 & \qquad \ga_2 \ga_3^2 & =2 &\qquad \ga_3 \gamma^2 &=0 \\
\go \gamma \ga_i &=0&\qquad \ga_1 \ga_2 \ga_3 &=0&\ga_1^2 \ga_3 &= 0&\qquad \ga_2 \ga_3 \gamma &=0& \ga_3^2 \gamma &=0\\
\go \gamma^2 &=0& \ga_1 \ga_2 \gamma&=0 & \ga_1^2 \gamma &=0 & \ga_2 \gamma^2 & =2 & \ga_3^3&=0 \\
\go^2 \ga_i &=0& \ga_1 \ga_3^2 &= 0 &\ga_1^3 &=0 & \ga_2^2 \ga_3 &= 0 & \gamma^3&=0\\
\go^2 \gamma &=0& \ga_1 \ga_3 \gamma &= 2& & & \ga_2^2 \gamma &= 2& &\\
\go^3 &=2 & \ga_1 \gamma^2&=0 & & & \ga_2^3 &= 0& & \\
%& &\ga_1^2 \ga_2 &= 0& \ga_3 \gamma^2 &=0\\
%& &\ga_1^2 \ga_3 &= 0& \ga_3^2 \gamma &=0\\
%& &\ga_1^2 \gamma &=0& \ga_3^3&=0\\
%& &\ga_1^3 &=0 & \gamma^3&=0\\
\end{alignat*}
With these choices, $\go \gamma=0$ and
this is the only 2-cohomology class that pairs trivially with $\gamma$.

Now we let $M$ be a simply connected spin 6-manifold with $H^2(M)=H^2$,
cup product as described above, arbitrary $H^3$ and first Pontryagin
class $p_1 = 4 \go^2$%(this one does not work in general... have to
%fill out the question marks for $\ga_i^3$)
. Let $w, a_i, c_1$ be  a set of  closed forms representing $\go, \ga_i, \gamma$ and let $X$ be a circle bundle
over $M$ with connection form $\theta$ and first Chern
class $\gamma$ with $d\theta = c_1$. Then $X$ is spin, has first Pontryagin class $p_1 =
4\go^2$ and has degree 2 and 3 cohomology $H^2(X)=\langle \go, \ga_1,
\ga_2, \ga_3 \rangle$, $H^3(X) 
= H^3(M)\oplus \langle [\theta\tensor \go] \rangle$. The term $[\theta \tensor \go]$ is the one we are concerned about: As $\gamma\go =0$, in the form level we have $c_1\wedge w  = d\xi$, for some 3-form $\xi$ pulled back from $M$, hence $\theta \wedge w - \xi$ is a closed form. This form represents the cohomology class $[\theta \tensor \go]$.

Letting $\gf = -\theta \wedge w +\xi$ we see
that \gf\ induces a negative definite bilinear form on $H^2(X)$ as
$$\int_X \gf \wedge a_i\wedge a_j = \int_M - w \wedge a_i \wedge a_j = -2\delta_{ij},$$
and similarly
$$\int_X \gf \wedge w\wedge a_j = 0 \qquad\mbox{and}\qquad \int_X \gf \wedge w\wedge w =  -2.$$

Also,
$$\int_X \gf p_1 = \int_X -\theta w 4w^2 = -\int_M4w^3 = - 8.$$

Finally, since $\gamma$ pulls back to zero in $H^2(X)$, we can define the
Massey products $\langle \go, \ga_1, \ga_2 \rangle$, $\langle \ga_1,
\ga_2,\ga_3\rangle$, $\langle \ga_2, \ga_3, \go \rangle$ and $\langle \ga_3,
\go, \ga_1\rangle$. To prove that $X$ is not
formal we compute $\langle \ga_1, \ga_2,\ga_3 \rangle$:
$$ a_1 a_2 = da_{12},\qquad a_2 a_3 = d(\theta a_1 - a_{23}),$$
where $a_{ij}$ are pull-backs from the base. So
$$\int_X \langle \ga_1, \ga_2,\ga_3 \rangle \go = \int_X a_{12} a_3 w - a_1
(\theta a_1-a_{23})w = - \int_X \theta  w a_1^2 = - \int_M w a_1^2 = - 2.$$
Which means that, for these choices, the Massey product pairs
nontrivially with the cohomology class $\go$ and therefore is a closed
nonexact form. One can also check that different choices {\it
  keep the integral above unchanged} so the Massey product does not
vanish.
\end{ex}

\begin{remark}
If we deal with the $Spin(7)$ case in a similar fashion, i.e., stripping off the Riemannian structure and working only with the implied topological properties, circle bundles will not provide possible examples of nonformal $Spin(7)$-manifolds. This is because $Spin(7)$-manifolds have $\hat{A}$-genus 1 and by a result of Atiyah and Hirzebruch \cite{AH70} if a compact connected Lie group acts differentiably and non-trivially on a compact  orientable spin manifold $X$, then $\hat A(X)=0$.
\end{remark}

\bibliographystyle{abbrv}
\bibliography{references}
\end{document}